\documentclass[a4paper]{report}
\usepackage{etoolbox}
\usepackage{t1enc}
\usepackage{amsmath}
\usepackage{amsfonts}
\usepackage{mathrsfs}
\usepackage{graphicx}
\usepackage{psfrag}
\usepackage{subfigure}
\usepackage{setspace}
\usepackage{cite}
\usepackage{amssymb}
\usepackage[hidelinks]{hyperref}
\usepackage[usenames,dvipsnames]{pstricks}
\usepackage{epsfig}
\usepackage{epstopdf}
\usepackage{wrapfig}
\usepackage{multirow}
\usepackage{hyperref}
\usepackage{xcolor}
\usepackage[english]{babel}
\patchcmd{\thebibliography}{\chapter*}{\section*}{}{}
\begin{document}
    \begin{singlespace}
    \begin{center}
    \begin{spacing}{0.8}
        \Large{\textsc{The Fourier transform of the non-trivial zeros of the zeta function\\}}
        \vspace{0.2in}
        \large{Levente Csoka\\}
				\vspace{0.1in}
				levente.csoka@skk.nyme.hu\\
        \vspace{0.2in}
        University of Sopron, 9400 Sopron, Bajcsy Zs. e. u. 4. Hungary\\
        \line(1,0){300}
    \end{spacing}
    \end{center}
    \end{singlespace}
\noindent
\textsc{ABSTRACT.} \emph{The non-trivial zeros of the Riemann zeta function and the prime numbers can be plotted by a modified von Mangoldt function. The series of non-trivial zeta zeros and prime numbers can be given explicitly by superposition of harmonic waves. The Fourier transform of the modified von Mangoldt functions shows interesting connections between the series. The Hilbert-Pólya conjecture predicts that the Riemann hypothesis is true, because the zeros of the zeta function correspond to eigenvalues of a positive operator and this idea encouraged to investigate the eigenvalues itself in a series. The Fourier transform computations is verifying the Riemann hypothesis and give evidence for additional conjecture that those zeros and prime numbers arranged in series that lie in the critical $\frac{1}{2}$ positive upper half plane and over the positive integers, respectively.}\\

Recent years have seen an explosion of research stating connection between areas of physics and mathematics \cite{1}. Hardy and Littlewood gave the first proof that infinitely many of the zeros are on the $\frac{1}{2}$ line\cite{2}. Hilbert and Pólya independently suggested a Hermitian operator whose eigenvalues are the non-trivial zeros. of Eq. \ref{eq1}. Pólya investigated the Fourier transform of Eq. \ref{eq1} as a consequence of the zeros's location of a functional equations. Physicist working in probability found that zeta function arises as an expectation in a moment of a Brownian bridge\cite{3}. Here, a new series is created from the non-trivial zeta zeros using the modified von Mangoldt function and then apply the Fourier transform to see the connection between the zeta zero members. To describe the numerical result, we consider the non-trivial zeros of the zeta function, by denoting
\begin{flalign}
    &s=\frac{1}{2}+it,\ \text{($i$ is real)}.&
		\label{eq1}
\end{flalign}

Let us consider a special, simple case of von Mangoldt function at for $k=$ non-trivial zeros of zeta function, 
\begin{flalign}
&L(n) = \left\{
\begin{array}{l l}
\ln\left(e\right) & \quad \text{if $n=$ non-trivial zeros},\\
0 & \quad \text{otherwise}.
\end{array} \right.&
\label{eq2}
\end{flalign}
The function $L(n)$ is a representative way to introduce the set of non-trivial zeta zeros with weight of unity attached to the location of a non-trivial zeta zeros. 
The discrete Fourier transform of the modified von Mangoldt function $L(n)$ gives a spectrum with periodic, real parts as spikes at the frequency axis ordinates. Hence the modified von Mangoldt function $L(n)$ and the non-trivial zeta zeros series can be approximated by periodic sin or cos functions. Let us consider that natural numbers form a discrete function $g(n)\rightarrow g(n_k)$ and $g_k\equiv g(n_k)$ by uniform sampling at every sampling point $n_k\equiv k\Delta$, where $\Delta$ is the sampling period and $k=2,3,..,L-1$. If the subset of discrete non-trivial zeta zeros numbers is generated by the uniform sampling of a wider discrete set, it is possible to obtain the Discrete Fourier Transform (DFT) of $L(n)$. The sequence of $L(n)$ non-trivial zeta zeros function is transformed into another sequence of $N$ complex numbers according to the following equation: 
\begin{flalign}
    &F\left\{L(n)\right\}=X\left(\nu\right)&
		\label{eq3}
\end{flalign}
where $L(n)$ is the modified von Mangoldt function. The operator $F$ is defined as:
\begin{flalign}
    &X\left(\nu\right)=\sum\limits_{n=1}^{N-1}L\left(k\Delta\right)e^{-i\left(2\pi\nu\right)\left(k\Delta\right)},&
		\label{eq4}
\end{flalign}
where $\nu=l\Delta f$, $l=0,1,2,…,N-1$ and $\Delta f=\frac{1}{L}$, $L=$length of the modified von Mangoldt function.
\setcounter{chapter}{1}
The amplitude spectrum of Eq. \ref{eq4}, describes completely the discrete-time Fourier transform of the $L(n)$ non-trivial zeta zeros function of quasi periodic sequences, which comprises only discrete frequency component. The important properties of the amplitude spectrum are as follows:
\begin{enumerate}
\item the ratio of consecutive frequencies $(f_t,f_{t+1})$ and maximum frequency $(f_{max})$ converges $\lim_{t\to\frac{f_s}{2}}\frac{f_{t+1}}{f_t}=\frac{f_t+\frac{1}{L}}{f_t}=\frac{f_{max}}{f_t}=1$, where $f_s=\frac{1}{\Delta}$ represents the positive half-length of amplitude spectrum. This is similar to the prime number theorem: $\lim_{x\to\infty}\frac{\pi\left(x\right)\ln\left(x\right)}{x}=1$\cite{4};
\item The reciprocate of the frequencies $(f_t)$ converges $\lim_{t\to \frac{f_s}{2}}\frac{1}{f_t}=\frac{1}{f_s}$;
\item The consecutive frequencies on the DFT amplitude spectrum describe a parabolic or Fermat's spiral (Fig. \ref{fig:fig1}), starting from 0 with increasing modulus (polar form $r=af_t$), where $a=1$, independently from the length of the non-trivial zeta zeros function and its parametric equations can be written as: $\chi\left(f_t\right)=f_t\cos\left(f_t\right), y\left(f_t\right)=f_t\sin\left(f_t\right)$;
\begin{figure}[h]
\centering
\includegraphics[width=40mm,height=40mm]{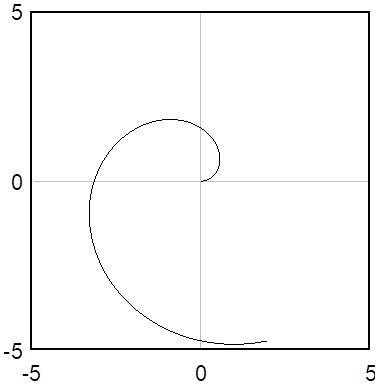}
\caption{Parabolic or Fermat spiral using the consecutive frequencies. The arc length of the spiral is a function of the sampling interval.}
\label{fig:fig1}
\end{figure}
\item The DFT spectrum of the non-trivial zeta zeros exhibits highly ordered and symmetrical frequency distribution (Fig. \ref{fig:fig3});
\item The DFT of the modified von Mangoldt function will be periodic (Fig. \ref{fig:fig2}) independently the length of integer sequences used:
\begin{flalign*}
&(\nu\geq N,\nu = 0, 1, ..., N-1)&
\end{flalign*}
\begin{flalign}
    &X\left(\nu\right)=X\left(\nu+N\right)=X\left(\nu+2N\right),...,=X\left(\nu+zN\right),&
\end{flalign}
\begin{flalign}
    &X\left(zN+\nu\right)=\sum\limits_{n=0}^{N-1}L\left(n\right)e^{-i\left(2\pi n\nu\right)\left(zN+\nu\right)}=\sum\limits_{n=0}^{N-1}L\left(n\right)e^{-i\left(2\pi n\nu\right)}e^{-i2\pi z n},&\\
		&z\text{ and }n\text{ are integers and }e^{-i2\pi z n}=1,&
\end{flalign}
\begin{flalign}
    &X\left(\nu+zN\right)=\sum\limits_{n=0}^{N-1}L\left(n\right)e^{-i\left(2\pi n\nu\right)}=X\left(\nu\right).&
\end{flalign}
\end{enumerate}
\section{Spectral evidence}
\begin{figure}[h]
\centering
\includegraphics[width=45mm,height=40mm]{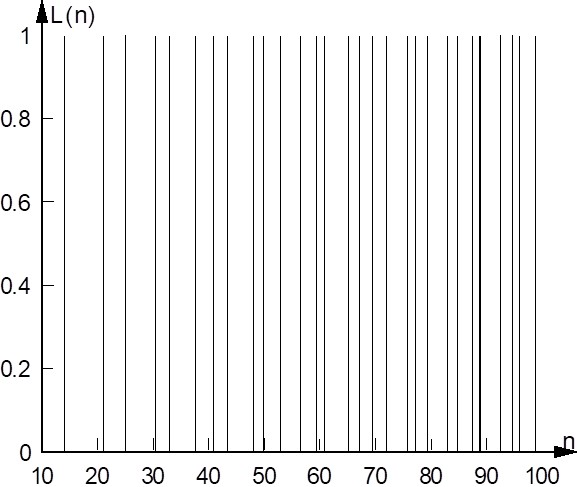}
\caption{Representation of the modified von Mangoldt function up to 100. Spikes appear at the location of non-trivial zeta zeros numbers.}
\label{fig:fig2}
\end{figure}

\begin{figure}[h]
\centering
\includegraphics[width=45mm,height=40mm]{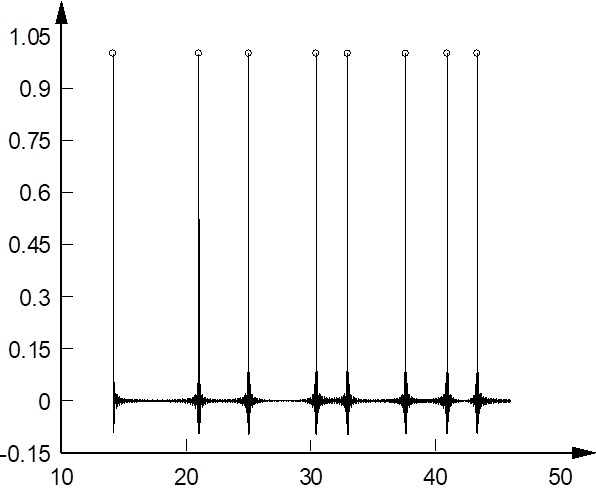}
\caption{Reconstruction of the modified von Mangoldt function from the Fourier decomposed periodic waves between 10-50. The circles at the top of the spikes indicates the location of non-trivial zeta zeros numbers.}
\label{fig:fig3}
\end{figure}

For completeness we give the compelling spectral evidence for the truth that Fourier decomposition of the modified von Mangoldt function is periodic. Fourier transformation has not been used for studying the behaviour of non-trivial zeta zeros numbers and it was employed here to reveal the periodic distribution of the non-trivial zeta zeros. Figure \ref{fig:fig2} shows that the corresponding amplitude spectrum of the modified von Mangoldt function exhibits a well-defined periodic signature. This suggests that the series of non-trivial zeta zeros numbers is not an arbitrary function. Creating the superposition from the Fermat spiral frequencies by adding their amplitude and phase together, the original non-trivial zeta zeros function can be restored. 
Creating the superposition from the Fermat spiral frequencies by adding their amplitude and phase together, the original non-trivial zeta zeros function can be arising.
\begin{flalign}
    &L\left(n\right)=\sum\limits_{f_t=0}^{f_{max}}A_t\sin\left(f_t+\omega n\right)&
\end{flalign}
From the wave sequence above, it can be seen that all non-trivial zeta zeros able to show some arithmetical progression, pattern of non-trivial zeta zeros numbers. Finding long arithmetic progressions in the non-trivial zeros attracted interest in the last decades as well. We can identify certain peaks on Fig. \ref{fig:fig2},which are the reciprocate of distances of the non-trivial zeta zeros in the arithmetic progressions.  
\section{Conclusion}
We have presented a novel method for the detection of non-trivial zeta zeros repetitions on the set of natural numbers. The method is based on a DFT of the modified von Mangoldt function. The relevance of the method for determination of the conformational sequences is described. The sequence of the non-trivial zeta zeros seems correlating with the unconditional theory for the spacing correlations of characteristic roots of zeta function developed by Katz and Sarnak\cite{5}. 
\renewcommand{\bibname}{References}

\end{document}